\RequirePackage{fix-cm}
\documentclass[smallextended]{svjour3}       
\smartqed  
\usepackage{graphicx}
\usepackage{amsmath}
\usepackage{amssymb}
\usepackage{float}
\usepackage{adjustbox}
\begin{document}

\title{Pythagorean Triples, B\'{e}zout Coefficients and the Erd\H{o}s-Straus Conjecture 
}


\author{Kyle Bradford
}


\institute{K. Bradford \at
              1332 Southern Dr, Statesboro, GA 30458  \\
              Tel.: +1-703-314-8920\\
              \email{kbradford@georgiasouthern.edu}           
}

\date{Received: date / Accepted: date}

\maketitle

\begin{abstract}
\ 

\keywords{}
\end{abstract}

\section{Introduction} \label{sec: intro}

The Erd\H{o}s-Straus conjecture was introduced in 1948 by Paul Erd\H{o}s and Ernst Straus.  It is one of many conjectures that involve Egyptian fractions.  These are sums of positive rational numbers with unit numerator.  As outlined in  \cite{Bra1}, to prove the Erd\H{o}s-Straus conjecture, it suffices to show that for each prime $p$ there exist positive integers $x \leq y \leq z$  so that the following diophantine equation is satisfied: 

\begin{equation} \label{eq: one}
\frac{4}{p} = \frac{1}{x} + \frac{1}{y} + \frac{1}{z}.
\end{equation}

\ 

\noindent If a solution exists for a given prime $p$, then by insisting that $x \leq y \leq z$  it was shown in \cite{Bra1} that $p$  cannot divide $x$, $p$ must divide $z$  and $p$  sometimes divides $y$.  This necessarily means that for every prime $p$, it is impossible to have $x=y=z$  simultaneously.  Also all solutions to (\ref{eq: one})  with either $x=y$  or $y=z$  are of the following form:

\begin{alignat*}{4}
\bullet &\quad p=2  &&\qquad x =1 &&\qquad y=2& &\qquad z=2 \\
\bullet &\quad p \equiv 3 \text{ mod } 4 &&\qquad x = \frac{p+1}{2} &&\qquad y = \frac{p+1}{2} &&\qquad z = \frac{p(p+1)}{4} \\
 \bullet &\quad p \equiv 3 \text{ mod } 4 &&\qquad x = \frac{p+1}{4} &&\qquad y = \frac{p(p+1)}{2} &&\qquad z = \frac{p(p+1)}{2}. \\
\end{alignat*}

\

\noindent  I will call these trivial solutions to (\ref{eq: one})  and focus on finding non-trivial solutions where $x<y<z$.  For each odd prime $p$  it is unknown whether or not a non-trivial solution to (\ref{eq: one}) exists.  If one exists, it was shown in \cite{Bra1} that

\begin{equation} \label{eq: two} 
z = \frac{xyp}{\gcd(p,y) \gcd(xy, x+y)}.
\end{equation}

\ 

\noindent Similarly it can be shown for an odd prime $p$  that if a non-trivial solution to (\ref{eq: one}) exists, then

\begin{equation} \label{eq: three} 
y = \frac{xz \gcd(p,y)}{p \gcd(xz, x+z)}
\end{equation}

\ 

\noindent and

\begin{equation} \label{eq: four} 
x = \frac{yz}{p \gcd(yz, y+z)}.
\end{equation}

\ 

\noindent Manipulating (\ref{eq: one}), (\ref{eq: two}), (\ref{eq: three}) and (\ref{eq: four}) leads to necessary equations for solutions to the Erd\H{o}s-Straus conjecture.  For a given odd prime $p$  and positive integers $x < y$, it is necessary that 

\begin{equation} \label{eq: five}
2 \left( \frac{2xy}{\gcd(xy, x+y)} \right) - \frac{p}{y-x} \left( \frac{y^{2} - x^{2}}{\gcd(xy, x+y)} \right) = \gcd(p,y),
\end{equation}

\ 

\noindent for a given odd prime $p$  and positive integers $x < z$, it is necessary that 

\begin{equation} \label{eq: six}
2 \left( \frac{2xz}{\gcd(xz, x+z)} \right) - \frac{p}{z-x} \left( \frac{z^{2} - x^{2}}{\gcd(xz, x+z)} \right) = \frac{p^{2}}{\gcd(p,y)},
\end{equation}

\ 

\noindent and for a given odd prime $p$  and positive integers $y < z$, it is necessary that

\begin{equation} \label{eq: seven}
2 \left( \frac{2yz}{\gcd(yz, y+z)} \right) - \frac{p}{z-y} \left( \frac{z^{2} - y^{2}}{\gcd(yz, y+z)} \right) = p^{2}.
\end{equation}

\ 

\noindent The reason to focus on the expressions in (\ref{eq: five}), (\ref{eq: six}) and (\ref{eq: seven}) is to highlight a connection to the work in \cite{Gue1}.  For example, by letting

\begin{equation} \label{eq: twelve}
\begin{split}
A &= \frac{2xy}{\gcd(xy, x+y)} \\
B &= \frac{y^{2} - x^{2}}{\gcd(xy, x+y)} 
\end{split}
\end{equation}

\ 

\noindent we have the two smaller legs of a Pythagorean triple.  I will outline these observations in greater detail in the following section.

\section{Pythagorean Triples} \label{sec: main}

\noindent A non-trivial Pythagorean triple $(A,B,C)$  has $A,B,C \in \mathbb{Z}^{+}$  so that $A^{2}+B^{2}=C^{2}$.  Either all three terms are even or only one of $A$ or $B$  is even, so let $A$  be an even term.  Note that this does not necessarily imply that $A<B$.  By insisting that $A$  is even, some unique examples Pythagorean triples under my insistence are $(4,3,5), (12,5,13), (6,8,10)$  and $(8,6,10)$.  Here I am treating $(6,8,10)$  and $(8,6,10)$  as different Pythagorean triples. \\

\begin{figure}[H] \label{fig: one}
\begin{center}
\fbox{\includegraphics[scale=0.40]{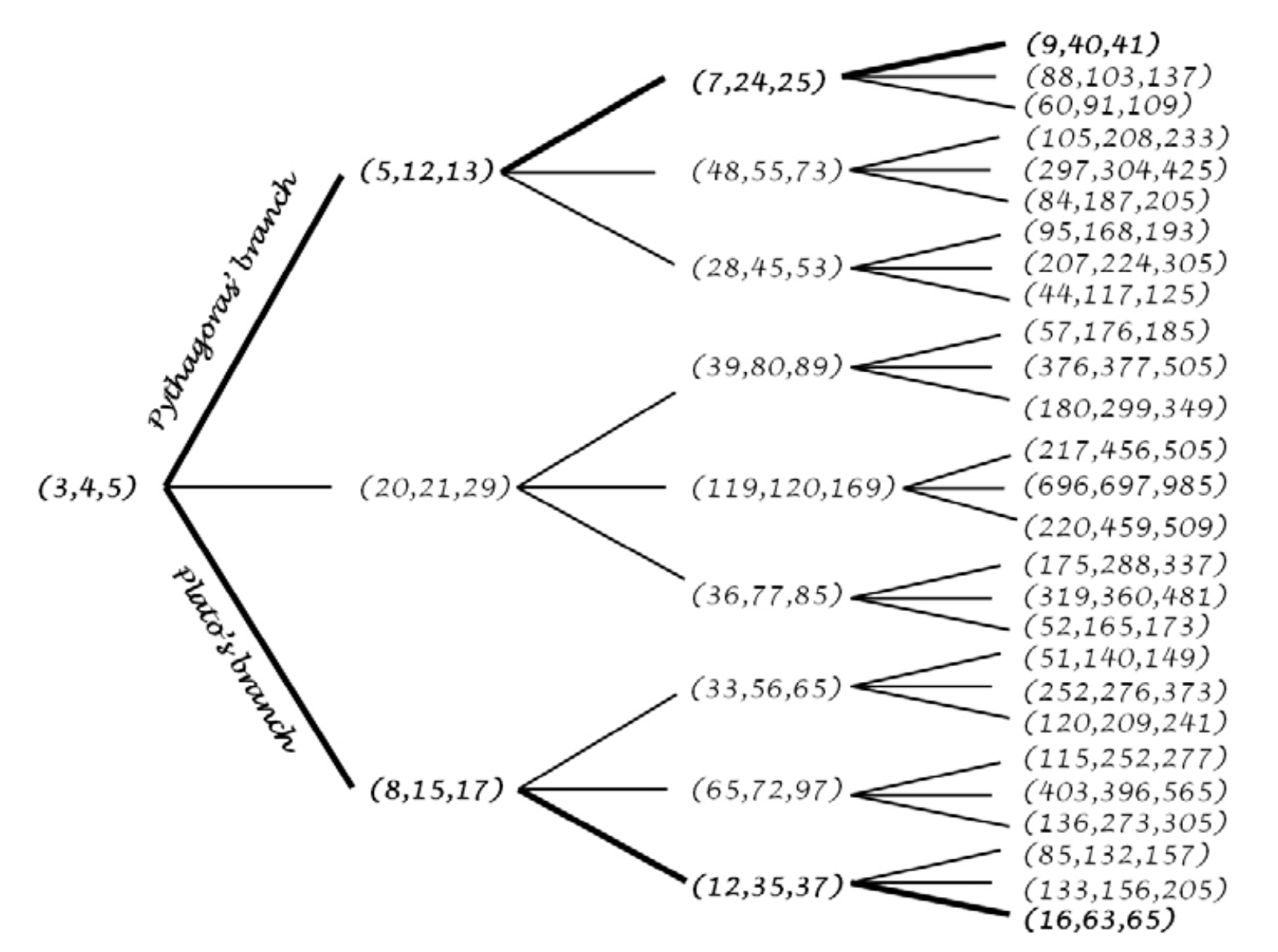}}
\end{center}
\caption{This image of a Berggren Tree is thanks to Luis Teia \cite{Tei1}.}
\end{figure}

\noindent A primitive Pythagorean triple has $A$  and $B$  coprime.  All primitive Pythagorean triples can be connected through a Berggren Tree \cite{Ber1}.  Figure \ref{fig: one} shows the history of the development of Pythagorean triples.  In this tree they do not insist that $A$  is even, rather that $A<B<C$.  For example, you can see that Pythagoras' branch has triples with C=B+1, Plato's branch has $C=B+2$, and the middle branch has $B=A+1$. \\

\noindent This tree has a root that can be expressed as a column vector $(3,4,5)^{T}$  and you can find the coefficients of any other Pythagorean triple through repeated left-hand multiplication of combinations of the following matrices:

$$
\begin{bmatrix}
1 & -2 & 2 \\
2 & -1 & 2 \\ 
2 & -2 & 3
\end{bmatrix} \qquad 
\begin{bmatrix}
1 & 2 & 2 \\
2 & 1 & 2 \\ 
2 & 2 & 3
\end{bmatrix} \qquad
\begin{bmatrix}
-1 & 2 & 2 \\
-2 & 1 & 2 \\ 
-2 & 2 & 3
\end{bmatrix}.
$$

\noindent One goal is to find a relationship between the odd primes and specific triples in this tree.  If a pattern can be found, then a solution to the Erd\H{o}s-Straus conjecture can be found.   

\subsection{\underline{Pythagorean triples of the first type}}

\noindent Given an odd prime $p$  and a non-trivial solution to (\ref{eq: one}) the formulas in (\ref{eq: twelve}) provide the smaller legs of a Pythagorean triple.  These will be called Pythagorean triples of the first type.  Figure \ref{fig: two}  gives these values for the non-trivial solutions for primes less than $19$.  These triples are not necessarily primitive.  In fact it can be shown that

$$  \gcd(A,B) = \gcd(x,y,z) \gcd \left( 2, \frac{y^{2} - x^{2}}{\gcd^{2}(x,y)} \right).$$

\ 

\noindent Manipulating (\ref{eq: five})  and (\ref{eq: twelve})  leads to the following expressions

\begin{equation} \label{eq: thirteen}
\begin{split}
x &= \frac{A-B+C}{2 \left( \frac{(2A- \gcd(p,y))}{p} \right)} \\
y &= \frac{A+B+C}{2 \left( \frac{(2A-\gcd(p,y))}{p} \right)} \\
z &= \frac{Ap}{2 \gcd(y,p)}.
\end{split}
\end{equation}

\ 

\noindent The equations in (\ref{eq: thirteen}) help us in the following way:  given an odd prime $p$  and a Pythagorean triple $(A,B,C)$, depending on the whether or not we are seeking solutions with $p$  dividing $y$, we have two sets of possible equations to find non-trivial solutions to (\ref{eq: one}).

\subsection{\underline{Pythagorean triples of the second type}}

\noindent Letting

\begin{equation} \label{eq: fourteen}
\begin{split}
A &= \frac{2xz}{\gcd(xz, x+z)} \\
B &= \frac{z^{2} - x^{2}}{\gcd(xz, x+z)} 
\end{split}
\end{equation}

\ 

\noindent we also have the two smaller legs of a Pythagorean triple. These will be called Pythagorean triples of the second type.  Figure \ref{fig: three}  gives these values for the non-trivial solutions for primes less than $19$.  These triples are also not necessarily primitive.  In fact it can be shown that

$$  \gcd(A,B) = \gcd(x,y,z) \gcd \left( 2, \frac{z^{2} - x^{2}}{\gcd^{2}(x,z)} \right).$$

\ 

\noindent Manipulating (\ref{eq: six})  and (\ref{eq: fourteen})  leads to the following expressions

\begin{equation} \label{eq: seventeen}
\begin{split}
x &= \frac{A-B+C}{2 \left( \frac{\left(2A- \frac{p^{2}}{\gcd(p,y)} \right)}{p} \right)} \\
y &= \frac{A\gcd(y,p)}{2 p} \\
z &= \frac{A+B+C}{2 \left( \frac{\left(2A-\frac{p^{2}}{\gcd(p,y)} \right)}{p} \right)}.
\end{split}
\end{equation}

\ 

\noindent The equations in (\ref{eq: seventeen}) help us in the following way:  given an odd prime $p$  and a Pythagorean triple $(A,B,C)$, depending on the whether or not we are seeking solutions with $p$  dividing $y$, we have two sets of possible equations to find non-trivial solutions to (\ref{eq: one}).

\subsection{\underline{Pythagorean triples of the third type}}

Letting

\begin{equation} \label{eq: fifteen}
\begin{split}
A &= \frac{2yz}{\gcd(yz, y+z)} \\
B &= \frac{z^{2} - y^{2}}{\gcd(yz, y+z)} 
\end{split}
\end{equation}

\ 

\noindent we also have the two smaller legs of a Pythagorean triple. These will be called Pythagorean triples of the third type.  Figure \ref{fig: four}  gives these values for the non-trivial solutions for primes less than $19$.  These triples are also not necessarily primitive.  In fact it can be shown that

$$  \gcd(A,B) = \gcd(p,y) \gcd(x,y,z) \gcd \left( 2, \frac{z^{2} - y^{2}}{\gcd^{2}(y,z)} \right).$$

\ 

\noindent Manipulating (\ref{eq: seven})  and (\ref{eq: fifteen})  leads to the following expressions

\begin{equation} \label{eq: sixteen}
\begin{split}
x &= \frac{A}{2p} \\
y &= \frac{A-B+C}{2 \left( \frac{(2A- p^{2})}{p} \right)} \\
z &= \frac{A+B+C}{2 \left( \frac{(2A-p^{2})}{p} \right)}.
\end{split}
\end{equation}

\ 

\noindent The equations in (\ref{eq: sixteen}) help us in the following way:  given an odd prime $p$  and a Pythagorean triple $(A,B,C)$, we have a set of equations to find non-trivial solutions to (\ref{eq: one}).

\section{B\'{e}zout Coefficients} \label{sec: Bezout}

\noindent Given a prime $p$  a solution to the conjecture requires finding positive integers $x<y$  that satisfy (\ref{eq: five}).  Making the substitutions allows us to find a Pythagorean triple smaller legs $A$  and $B$  instead.  There is one more fundamental way to decompose this problem and it is discussed in this section.  Ultimately we will be looking for two other integer values instead.  While the complexity is not reduced, it gives insight to the solution values as roots of polynomials. \\

\noindent Let $p$ be an odd prime.  Solutions to (\ref{eq: one}) with $\gcd(p,y)=1$  were called type I solutions in \cite{Bra1}.  There are necessary conditions for these types of solutions.  To further understand these necessary conditions, the first goal is to find an expression for all values $b,c \in \mathbb{N}$  such that 

\begin{equation} \label{eq: nine}
4c - pb = \gcd( c,b).
\end{equation}  

\ 

\noindent Notice that if we can find $x,y \in \mathbb{N}$  so that $b=x+y$  and $c=xy$, then we have type I solutions to (\ref{eq: three}). \\

\noindent Let $b_{1} = 4 \left\lceil p \slash 4 \right\rceil - p$  and let $c_{1} = (b_{1} p+1) \slash 4$. \\

\noindent Notice that both $b_{1}$  and $c_{1}$  are positive integers and $$ 4 c_{1} + (-p)b_{1} = 1.$$

\ 

\noindent This by definition shows that $b_{1}$  and $c_{1}$  are coprime.  Indeed we see that  $$4 c_{1} - p b_{1} = \gcd( c_{1}, b_{1})$$

\ 

\noindent where $| b_{1}| \leq 4$  and $|c_{1}| \leq p$.   This guarantees that these are the minimal B\'{e}zout coefficients. \\

\noindent  For $k \in \mathbb{N}$  let $b_{k} = 4(k-1) + b_{1}$  and let $c_{k} = p(k-1) + c_{1}$.  \\

\noindent Notice that $b_{k}$  and $c_{k} = (b_{k}p + 1) \slash 4$  are always positive integers.  Also notice that $b_{k}$  and $c_{k}$  are by definition coprime for each $k \in \mathbb{N}$. \\

\noindent It should also be clear that these are the only values so that $$4 c_{k} - p b_{k} = \gcd( c_{k}, b_{k}).$$

\ 

\noindent B\'{e}zout's identity tells us that multiples of these coefficients are the only solutions to (\ref{eq: four}).  In other words the solution set for possible type I solutions have $x,y \in \mathbb{N}$  so that  $$(x+y, xy) \in \{ (mb_{k}, mc_{k}) | m,k \in \mathbb{N}, \gcd(b_{k},c_{k}) =1 \}. $$

\ 

\noindent Now let p be an any prime. Solutions to (\ref{eq: one}) with $\gcd(p,y)=p$  were called type II solutions in \cite{Bra1}.  There are also necessary conditions for these types of solutions.  To further understand these necessary conditions, the next goal is to find an expression for all values $b,c \in \mathbb{N}$  such that 

\begin{equation} \label{eq: ten}
4c - pb = p \gcd( c,b).
\end{equation}  

\ 

\noindent Similarly notice that if we can find $x,y \in \mathbb{N}$  so that $b=x+y$  and $c=xy$, then we have type II solutions to (\ref{eq: three}). \\

\noindent Let $b_{1} = 3$   and let $c_{1} = p$. \\

\noindent  For $k \in \mathbb{N}$  let $b_{k} = 4(k-1)+b_{1}$  and let $c_{k} = p(k-1) + c_{1}$.  \\

\noindent Notice that $b_{k}$  and $c_{k} = (b_{k} + 1)p \slash 4$  are positive integers for all $k \in \mathbb{N}$.  Also notice that $b_{k}$  and $c_{k}$  are coprime when $p \neq 3$.  If $p = 3$,  then $b_{k}$  and $c_{k}$  are coprime for $k \not\equiv 1 \mod 3$.  \\

\noindent It should be clear that these are the only values so that $$4 c_{k} - p b_{k} = p \gcd( c_{k}, b_{k}).$$

\ 

\noindent Using a similar reasoning we can conclude that multiples of these coefficients are the only solutions to (\ref{eq: five}). The solution set for possible type II solutions have $x,y \in \mathbb{N}$  so that $$ (x+y, xy) \in \{ (mb_{k}, mc_{k}) | m,k \in \mathbb{N}, \gcd(b_{k},c_{k}) =1 \}. $$

\ 

\noindent In general, to solve the Erd\H{o}s-Straus conjecture we are looking for solutions like the ones outlined above where $x+y = mb_{k}$  and $xy = mc_{k}$  for some $m,k \in \mathbb{N}$  with $\gcd(b_{k},c_{k})=1$.  As these two equations are symmetric with respect to $x$  and $y$, we see that these two values are going to be the roots of a polynomial $T^2 - mb_{k}T +mc_{k}$.  For each prime $p$, if (\ref{eq: one}). \\

\noindent The roots of this polynomial are going to be integers if and only if $$m^{2} b_{k}^{2} - 4mc_{k}$$  

\ 

\noindent is a square number. \\

\noindent For type I solutions this reduces to finding $m,k \in \mathbb{N}$  so that $$ m^{2} (4k-(4-4 \left\lceil p \slash 4 \right\rceil + p))^{2} - pm(4k-(4-4 \left\lceil p \slash 4 \right\rceil + p)) -m$$

\ 

\noindent is a square number. \\

\noindent For type II solutions this reduces to finding $m,k \in \mathbb{N}$  so that $$ m^{2} (4k-1)^{2} - 4pmk$$

\ 

\noindent is a square number, where $k \not\equiv 1 \mod 3$ for $p =3$. \\

\noindent Consider a type II solution.  Let $y^{*} = y \slash p$  and $z^{*} = z \slash p$. \\

\noindent One can show that $$  4y^{*}z^{*} - (y^{*}+z^{*}) = p \gcd(y^{*}z^{*}, y^{*}+z^{*}).$$

\ 

\noindent  Let $A^{*} = y^{*}z^{*} \slash \gcd(y^{*}z^{*}, y^{*}+z^{*})$  and let $B^{*} = (y^{*}+z^{*}) \slash  \gcd(y^{*}z^{*}, y^{*}+z^{*})$.  \\

\noindent Notice that $B^{*} = 4A^{*} - p$.  Because $p = 4\left( \left\lceil p \slash 4 \right\rceil \right) - \left( 4 \left\lceil p \slash 4 \right\rceil - p \right)$, we see that $B^{*} = 4\left(A^{*} - \left\lceil p \slash 4 \right\rceil \right) + \left( 4 \left\lceil p \slash 4 \right\rceil - p \right)$. Leting $m=  \gcd(y^{*}z^{*}, y^{*}+z^{*})$ and letting $A = (k-1) + \left\lceil p \slash 4 \right\rceil$,  if we can similarly show that $$m^{2} \left(4k - \left( 4 - 4 \left\lceil p \slash 4 \right\rceil + p \right) \right)^{2} - m\left( 4k - \left(4 - 4 \left\lceil p \slash 4 \right\rceil + p \right) \right)- mp$$

\ 

\noindent is a square number, then we would be able to find our solution.  As $A^{*}$ is an expression for $x$  for type II solutions, we see that $(k-1) \leq \left\lceil p \slash 4 \right\rceil$.

\noindent There are many other second degree polynomials in $m$  and $k$  that can be found and the ultimate goal will be to show that the expressions are square.  In \cite{Bra1} a functional relationship between $x$  and $y$  was found for type I solutions, but no such relationship has been found between the legs of the Pythagorean triple smaller legs or the values of $m$  and $k$.  The point of considering these different techniques is to develop a deeper understanding of the problem and to find creative ways to approach the problem. Hopefully this expository paper will bring more attention to the problem and the work in \cite{Gue1}.














\section{Appendix} \label{sec: appendix}

\begin{figure}[h] \label{fig: two}
\begin{adjustbox}{width=0.75 \textwidth,center}
\fbox{\begin{tabular}{c|ccc|ccc}
p & x & y & z & A & B & C \\ \hline
3 & 1 & 4 & 12 & 8 & 15 & 17 \\ \hline
5 & 2 & 4 & 20 & 8 & 6 & 10 \\ \hline
5 & 2 & 5 & 10 & 20 & 21 & 29 \\ \hline
7 & 2 & 15 & 210 & 60 & 221 & 229 \\ \hline
7 & 2 & 16 & 112 & 32 & 126 & 130 \\ \hline
7 & 2 & 18 & 63 & 18 & 80 & 82 \\ \hline
7 & 2 & 21 & 42 & 84 & 437 & 445 \\ \hline
7 & 3 & 6 & 14 & 4 & 3 & 5 \\ \hline
11 & 3 & 34 & 1122 & 204 & 1147 & 1165 \\ \hline
11 & 3 & 36 & 396 & 72 & 429 & 435 \\ \hline
11 & 3 & 42 & 154 & 28 & 195 & 197 \\ \hline
11 & 3 & 44 & 132 & 264 & 1927 & 1945 \\ \hline
11 & 4 & 9 & 396 & 72 & 65 & 97 \\ \hline
11 & 4 & 11 & 44 & 88 & 105 & 137 \\ \hline
11 & 4 & 12 & 33 & 6 & 8 & 10 \\ \hline
13 & 4 & 18 & 468 & 72 & 154 & 170 \\ \hline
13 & 4 & 20 & 130 & 20 & 48 & 52 \\ \hline
13 & 4 & 26 & 52 & 104 & 330 & 346 \\ \hline
13 & 5 & 10 & 130 & 20 & 15 & 25 \\ \hline
17 & 5 & 30 & 510 & 60 & 175 & 185 \\ \hline
17 & 5 & 34 & 170 & 340 & 1131 & 1181 \\ \hline
17 & 6 & 15 & 510 & 60 & 63 & 87 \\ \hline
17 & 6 & 17 & 102 & 204 & 253 & 325
\end{tabular}}
\end{adjustbox}
\caption{This table shows related Pythagorean triples of the first type for some non-trivial solutions to the Erd\H{o}s-Straus equation.}
\end{figure}

\begin{figure}[h] \label{fig: three}
\begin{adjustbox}{width=0.75 \textwidth, center}
\fbox{\begin{tabular}{c|ccc|ccc}
p & x & y & z & A & B & C \\ \hline
3 & 1 & 4 & 12 & 24 & 143 & 145 \\ \hline
5 & 2 & 4 & 20 & 40 & 198 & 202 \\ \hline
5 & 2 & 5 & 10 & 10 & 24 & 26 \\ \hline
7 & 2 & 15 & 210 & 210 & 11024 & 11026 \\ \hline
7 & 2 & 16 & 112 & 224 & 6270 & 6274 \\ \hline
7 & 2 & 18 & 63 & 252 & 3965 & 3973 \\ \hline
7 & 2 & 21 & 42 & 42 & 440 & 442 \\ \hline
7 & 3 & 6 & 14 & 84 & 187 & 205 \\ \hline
11 & 3 & 34 & 1122 & 748 & 139875 & 139877 \\ \hline
11 & 3 & 36 & 396 & 792 & 52269 & 52275 \\ \hline
11 & 3 & 42 & 154 & 924 & 23707 & 23725 \\ \hline
11 & 3 & 44 & 132 & 88 & 1935 & 1937 \\ \hline
11 & 4 & 9 & 396 & 198 & 9800 & 9802 \\ \hline
11 & 4 & 11 & 44 & 22 & 120 & 122 \\ \hline
11 & 4 & 12 & 33 & 264 & 1073 & 1105 \\ \hline
13 & 4 & 18 & 468 & 468 & 27376 & 27380 \\ \hline
13 & 4 & 20 & 130 & 520 & 8442 & 8458 \\ \hline
13 & 4 & 26 & 52 & 52 & 336 & 340 \\ \hline
13 & 5 & 10 & 130 & 260 & 3375 & 3385 \\ \hline
17 & 5 & 30 & 510 & 1020 & 52015 & 52025 \\ \hline
17 & 5 & 34 & 170 & 68 & 1155 & 1157 \\ \hline
17 & 6 & 15 & 510 & 510 & 21672 & 21678 \\ \hline
17 & 6 & 17 & 102 & 34 & 288 & 290
\end{tabular}}
\end{adjustbox}
\caption{This table shows related Pythagorean triples of the second type for some non-trivial solutions to the Erd\H{o}s-Straus equation.}
\end{figure}

\begin{figure}[h] \label{fig: four}
\begin{adjustbox}{width=0.75 \textwidth, center}
\fbox{\begin{tabular}{c|ccc|ccc}
p & x & y & z & A & B & C \\ \hline
3 & 1 & 4 & 12 & 6 & 8 & 10 \\ \hline
5 & 2 & 4 & 20 & 20 & 48 & 52 \\ \hline
5 & 2 & 5 & 10 & 20 & 15 & 25 \\ \hline
7 & 2 & 15 & 210 & 28 & 195 & 197 \\ \hline
7 & 2 & 16 & 112 & 28 & 96 & 100 \\ \hline
7 & 2 & 18 & 63 & 28 & 45 & 53 \\ \hline
7 & 2 & 21 & 42 & 28 & 21 & 35 \\ \hline
7 & 3 & 6 & 14 & 42 & 40 & 58 \\ \hline
11 & 3 & 34 & 1122 & 66 & 1088 & 1090 \\ \hline
11 & 3 & 36 & 396 & 66 & 360 & 366 \\ \hline
11 & 3 & 42 & 154 & 66 & 112 & 130 \\ \hline
11 & 3 & 44 & 132 & 66 & 88 & 110 \\ \hline
11 & 4 & 9 & 396 & 88 & 1935 & 1937 \\ \hline
11 & 4 & 11 & 44 & 88 & 165 & 187 \\ \hline
11 & 4 & 12 & 33 & 88 & 105 & 137 \\ \hline
13 & 4 & 18 & 468 & 104 & 1350 & 1354 \\ \hline
13 & 4 & 20 & 130 & 104 & 330 & 346 \\ \hline
13 & 4 & 26 & 52 & 104 & 78 & 130 \\ \hline
13 & 5 & 10 & 130 & 130 & 840 & 850 \\ \hline
17 & 5 & 30 & 510 & 170 & 1440 & 1450 \\ \hline
17 & 5 & 34 & 170 & 170 & 408 & 442 \\ \hline
17 & 6 & 15 & 510 & 204 & 3465 & 3471 \\ \hline
17 & 6 & 17 & 102 & 204 & 595 & 629
\end{tabular}}
\end{adjustbox}
\caption{This table shows related Pythagorean triples of the third type for some non-trivial solutions to the Erd\H{o}s-Straus equation.}
\end{figure}

\clearpage


\end{document}